\documentclass[12pt]{amsart}

\newtheorem{theorem}{Theorem}[section]
\newtheorem{lemma}[theorem]{Lemma}

\usepackage[dvips]{graphics}

\theoremstyle{definition}

\theoremstyle{remark}
\newtheorem{remark}[theorem]{Remark}

\theoremstyle{proposition}
\newtheorem{proposition}[theorem]{Proposition}

\theoremstyle{conjecture}

\theoremstyle{corollary}
\newtheorem{corollary}[theorem]{Corollary}

\theoremstyle{problem}

\numberwithin{equation}{section}

\large\normalsize

\begin{document}

\title{Asymptotic energy of lattices}

\author{Weigen Yan}
\address{School of Sciences, Jimei University,
Xiamen 361021, China} \email{weigenyan@263.net}
\thanks{The first author was supported in part by NSFC
Grant(10771086) and by Program for New Century Excellent Talents in
Fujian Province University.}

\author{Zuhe Zhang}
\address{Corresponding author, Department of Mathematics,
West Virginia University, Morgantown, WV 26506-6310, USA}
\email{zhzhang@math.wvu.edu}



\begin{abstract}
The energy of a simple graph $G$ arising in chemical physics,
denoted by $\mathcal E(G)$, is defined as the sum of the absolute
values of eigenvalues of $G$. As the dimer problem and spanning
trees problem in statistical physics, in this paper we propose the
energy per vertex problem for lattice systems. In general for a type
of lattices in statistical physics, to compute the entropy constant
with toroidal, cylindrical, Mobius-band, Klein-bottle, and free
boundary conditions are different tasks with different hardness and
may have different solution. We show that the energy per vertex of
plane lattices is independent on the toroidal, cylindrical,
Mobius-band, Klein-bottle, and free boundary conditions.
Particularly, the asymptotic formulae of energies of the triangular,
$3^3.4^2$, and hexagonal lattices with toroidal, cylindrical,
Mobius-band, Klein-bottle, and free boundary conditions are obtained
explicitly.
\end{abstract}
\maketitle
\section{Introduction}
Throughout this paper, we suppose that $G=(V(G),E(G))$ is a simple
graph with the vertex set $V(G)=\{v_1,v_2,\cdots,v_n\}$ and the edge
set $E(G)$, if not specified. The adjacency matrix of $G$ with $n$
vertices, denoted by $A(G)=(a_{ij})_{n\times n}$, is an $n\times n$
symmetric matrix such that $a_{ij}=1$ if vertices $v_i$ and $v_j$
are adjacent and $0$ otherwise. The characteristic polynomial of
$G$, denoted by $\phi(G,x)$, is defined as $\det(xI_n-A(G))$, where
$I_n$ is an identity matrix of order $n$. Denote the degree of
vertex $v_i$ of $G$ by $d_G(v_i)$. If $H$ is a subgraph of $G$, then
$G-E(H)$ denotes the subgraph obtained from $G$ by deleting all
edges in $H$.

Gutman \cite{Gutm78,Gutm99} defined the energy of a graph $G$ with
$n$ vertices, denoted by $\mathcal E(G)$, as $\mathcal
E(G)=\sum_{i=1}^{n}|\lambda_i(G)|,$ where $\lambda_i(G)$ are the
eigenvalues of the adjacency matrix of $G$.


In statistical physics the dimer problem considers the molecular
freedom (free energy per dimer) (see for example
\cite{Fisher61,Kast61,Kast63,Kast67,TF61,WFY06}), the spanning tree
problem considers the entropy of spanning trees
\cite{CW06,SW00,YZ08}, and the independent set problem considers the
entropy of independent sets \cite{CW98}. It is natural to consider
the chemical physics parameter - the energy of lattices from the
statistical physics¡¯ point of view. For the case of quadratic
lattices, to compute the energy per vertex is an easy task as
follows.

Suppose that $G_n, G_n'$, and $G_n''$ are the plane square lattices
with toroidal, cylindrical, and free boundary conditions,
respectively. That is, $G_n=C_n\times C_n, G_n'=P_n\times C_n$, and
$G_n''=P_n\times P_n$, where $C_n$ and $P_n$ denote the cycle and
the path with $n$ vertices, and $G\times H$ is the Cartesian product
of two graphs $G$ and $H$. Obviously, $\{G_n''\}$ is a sequence of
spanning subgraphs of the sequence $\{G_n'\}$ of finite graphs, and
$\{G_n'\}$ is a sequence of spanning subgraphs of the sequence
$\{G_n\}$ of finite graphs. Particularly,
$$\lim\limits_{n\rightarrow \infty} \frac{|\{v\in V(G_n''):
d_{G_n''}(v)=d_{G_n}(v)\}|}{|V(G_n)|}=\lim\limits_{n\rightarrow
\infty} \frac{|\{v\in V(G_n'):
d_{G_n'}(v)=d_{G_n}(v)\}|}{|V(G_n)|}=1,$$ that is, almost all
vertices of $G_n$ and $G_n''$ (resp. $G_n$ and $G_n'$) have the same
degrees. On the other hand, it is well known that the eigenvalues of
$G_n''$ (resp. $G_n'$ and $G_n$) are
$2\cos\frac{i\pi}{n+1}+2\cos\frac{j\pi}{n+1}, i,j=1,2\ldots,n$
(resp. $2\cos\frac{i\pi}{n+1}+2\cos\frac{2j\pi}{n}, i=1,2,\ldots,n,
j=0,1,2,\ldots,n-1$, and $2\cos\frac{2i\pi}{n}+2\cos\frac{2j\pi}{n},
i,j=0,1,\ldots,n-1$). Hence the energy per vertex of $G_n'', G_n'$,
and $G_n$ are defined as
$$\lim_{n\rightarrow \infty}\frac{\mathcal E(G_n'')}{|V(G_n'')|}=
\lim_{n\rightarrow \infty}\frac{\mathcal E(G_n'')}{n^2}=
\lim_{n\rightarrow
\infty}\frac{1}{n^2}\sum_{i,j=1}^{n}(|2\cos\frac{i\pi}{n+1}+2\cos\frac{j\pi}{n+1}|)
$$$$=2\int_{0}^1\int_{0}^1|\cos{\pi x}+\cos{\pi y}|d_xd_y
=\frac{2}{\pi^2}\int_{0}^{\pi}\int_{0}^{\pi}|\cos{
x}+\cos{y}|d_xd_y\approx 1.6211;$$
$$\lim_{n\rightarrow \infty}\frac{\mathcal E(G_n')}{|V(G_n')|}=
\lim_{n\rightarrow \infty}\frac{\mathcal E(G_n')}{n^2}=
\lim_{n\rightarrow
\infty}\frac{1}{n^2}\sum_{i=1}^{n}\sum_{j=0}^{n-1}(|2\cos\frac{i\pi}{n+1}+2\cos\frac{2j\pi}{n}|)
$$$$=2\int_{0}^1\int_{0}^1|\cos{\pi x}+\cos{2\pi y}|d_xd_y
=\frac{1}{\pi^2}\int_{0}^{\pi}\int_{0}^{2\pi}|\cos{
x}+\cos{y}|d_xd_y\approx 1.6211;$$
$$\lim_{n\rightarrow \infty}\frac{\mathcal E(G_n)}{|V(G_n)|}=
\lim_{n\rightarrow \infty}\frac{\mathcal E(G_n)}{n^2}=
\lim_{n\rightarrow
\infty}\frac{1}{n^2}\sum_{i,j=0}^{n-1}(|2\cos\frac{2i\pi}{n}+2\cos\frac{2j\pi}{n}|)
$$$$=2\int_{0}^1\int_{0}^1|\cos{2\pi x}+\cos{2\pi y}|d_xd_y
=\frac{1}{2\pi^2}\int_{0}^{2\pi}\int_{0}^{2\pi}|\cos{
x}+\cos{y}|d_xd_y\approx 1.6211,$$ implying $G_n, G_n'$, and $G_n''$
have the same asymptotic energy ($\approx 1.6211n^2$).

The phenomenon above is not accidental. In this paper we obtain the
asymptotic formulae of energies of triangular, $3^3.4^2$, and
hexagonal lattices with toroidal, cylindrical, Mobius-band,
Klein-bottle, and free boundary conditions. Our approach implies
that in general the energy per vertex of plane lattices is
independent of the boundary conditions.

\section{The asymptotic energies of some lattices}

\subsection{Graph asymptotic energy change due to edge deletion}

Let us record the following results firstly. Koolen and Moulton
\cite{Kool01} proved that if $2m\geq n$ and $G$ is a graph on $n$
vertices with $m$ edges, then
$$\mathcal E(G)\leq \frac{2m}{n}+\sqrt{(n-1)\left[2m-\left(\frac{2m}{n}\right)^2\right]},\eqno{(1)}$$
and if $2m\leq n$ and $G$ is a graph on $n$ vertices with $m$ edges,
then $$\mathcal E(G)\leq 2m.\eqno{(1')}$$

The following result is immediate from $(1)$ and $(1')$.
\begin{proposition}
Let $G$ be a graph with $m$ edges. Then
$$\mathcal E(G)\leq 2m.$$
\end{proposition}

Day and So \cite{DS07,DS08} first studied how the energy of a graph
changes when edges are removed. They found the following
\begin{lemma}[Day and So \cite{DS07,DS08}]
Let $H$ be an induced subgraph of a graph $G$. Then
$$\mathcal E(G)-\mathcal E(H)\leq \mathcal E(G-E(H))\leq \mathcal E(G)+\mathcal E(H).$$
\end{lemma}

With a similar method, we can prove the following:

\begin{lemma}
Let $H$ be a subgraph of a graph $G$. Then
$$|\mathcal E(G)-\mathcal E(H)|\leq \mathcal E(G-E(H))\leq \mathcal E(G)+\mathcal E(H).$$
\end{lemma}

Given two graphs $G$ and $H$ ($V(G)\cap V(H)$ may be disjoint),
Denoted by $\Delta(G,H)=|E(G)|+|E(H)|-2|E(G)\cap E(H)|$, i.e.,
$\Delta(G,H)$ equals the number of edges of symmetric difference of
$E(G)$ and $E(H)$.
\begin{theorem}
Suppose $\{G_n\}$ and $\{H_n\}$ are two sequences of graphs such
that $$\lim_{n\rightarrow \infty}\frac{\Delta(G_n,H_n)}{\mathcal
E(G_n)}=0.$$Then $$\lim_{n\rightarrow \infty}\frac{\mathcal
E(H_n)}{\mathcal E(G_n)}=1.$$
\end{theorem}
\begin{proof}
Let $F_n$ be the subgraph of $G_n$ or $H_n$ induced by $E(G_n)\cap
E(H_n)$. Note that
$$\left|\frac{\mathcal E(H_n)}{\mathcal
E(G_n)}-1\right|=\left|\frac{\mathcal E(H_n)-\mathcal
E(G_n)}{\mathcal E(G_n)}\right|=\left|\frac{\mathcal E(H_n)-\mathcal
E(F_n)+\mathcal E(F_n)-\mathcal E(G_n)}{\mathcal E(G_n)}\right|$$
$$\leq \left|\frac{\mathcal E(G_n)-\mathcal
E(F_n)}{\mathcal E(G_n)}\right|+\left|\frac{\mathcal E(H_n)-\mathcal
E(F_n)}{\mathcal E(G_n)}\right|.$$By Lemma 2.3 and Proposition 2.1,
$$|\mathcal E(G_n)-\mathcal E(F_n)|\leq \mathcal E(G_n-E(F_n))\leq
2|E(G_n)|-2|E(F_n)|,$$ $$|\mathcal E(H_n)-\mathcal E(F_n)|\leq
\mathcal E(H_n-E(F_n))\leq 2|E(H_n)|-2|E(F_n)|.$$ Hence
$$\left|\frac{\mathcal E(H_n)}{\mathcal
E(G_n)}-1\right|\leq \frac{2\Delta(G_n,H_n)}{\mathcal E(G_n)}$$
implying the theorem.
\end{proof}
\begin{corollary}
Suppose that $\{G_n\}$ is a sequence of finite simple graphs with
bounded average degree such that $\lim\limits_{n\rightarrow
\infty}|V(G_n)|=\infty$ and $\lim\limits_{n\rightarrow
\infty}\frac{\mathcal E(G_n)}{|V(G_n)|}=h\neq 0$. If $\{G_n'\}$ is a
sequence of spanning subgraphs of $\{G_n\}$ such that
$\lim\limits_{n\rightarrow \infty} \frac{|\{v\in V(G_n'):
d_{G_n'}(v)=d_{G_n}(v)\}|}{|V(G_n)|}=1$, then
$\lim\limits_{n\rightarrow \infty}\frac{\mathcal
E(G'_n)}{|V(G_n')|}=h$. That is, $G_n$ and $G_n'$ have the same
asymptotic energy.
\end{corollary}

A direct sequence of Corollary 2.5 is that $P_n\times P_n, P_n\times
C_n$, and $C_n\times C_n$ have the same asymptotic energy which is
shown in the introduction. More generally, by Corollary 2.5, we have
\begin{remark}
Suppose $G_i=P_n$ or $G_i=C_n$, $i=1,2,\ldots,k$, and $k$ is a
constant. If $n$ is sufficiently large, then the asymptotic energy
of the $n$-dimensional lattices
$$\mathcal E(G_1\times G_2\times...\times G_k)\approx\frac{2n^k}{\pi^k}
\int_{0}^{\pi}\int_{0}^{\pi}...\int_{0}^{\pi}\left(\sum_{i=1}^k|\cos
x_i|\right){\rm d}_{x_1}{\rm d}_{x_2}\ldots {\rm d}_{x_k}.$$
\end{remark}

\begin{remark}
Corollary 2.5 gives a method to calculate the asymptotic energy of a
graph with bounded average degree. Suppose that $\{G_n\}$ is a
sequence of finite simple graphs with bounded average degree. It is
difficult to calculate its asymptotic energy directly. We can find a
graph $G_n'$ with bounded average degree, which satisfies
$|V(G_n)|=|V(G_n')|$ and almost all vertices of $G_n$ and $G_n'$
have the same degrees. If we can compute the asymptotic energy of
$G_n'$ directly, then by Corollary 2.5, $G_n$ and $G_n'$ have the
same asymptotic energy. We will use this idea to calculate the
asymptotic energy of some graphs in the next subsections.
\end{remark}

\begin{figure}[htbp]
  \centering
 \scalebox{0.7}{\includegraphics{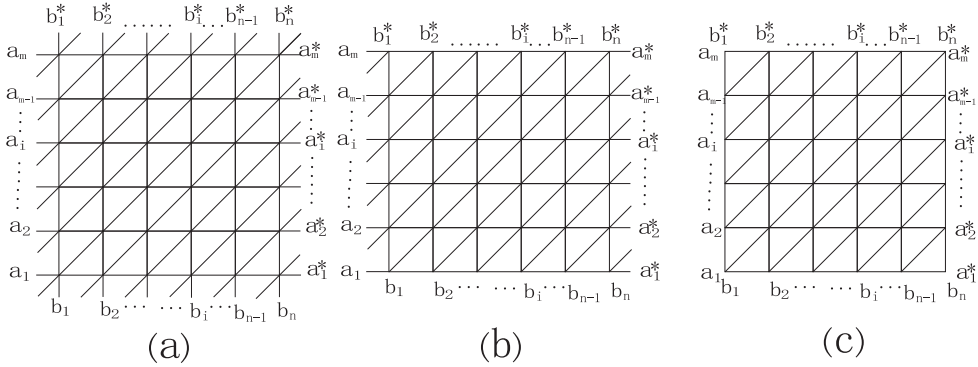}}
  \caption{\ (a)\ The triangular lattice $T^t(n,m)$ with toroidal boundary
  condition;
  \ (b)\ the triangular lattice $T^c(n,m)$ with cylindrical boundary condition;
  \ (c)\ the triangular lattice $T^f(n,m)$ with free boundary condition.}
\end{figure}

\subsection{The triangular lattice}

The triangular lattice with toroidal boundary condition, denoted by
$T^t(n,m)$, can be regarded as an $n\times m$ square lattice with
toroidal boundary condition with an additional diagonal edge added,
in the same way, to every square, see Figure 1(a), where
$(a_1,a_1^*),(a_2,a_2^*),\ldots,(a_m,a_m^*);(b_1,b_1^*),(b_2,b_2^*),\ldots,(b_n,b_n^*);
(b_1^*,b_2),(b_2^*,b_3),\\
\ldots,(b_{n-1}^*,b_n),(b_n^*,b_1)=(a_m^*,a_1);(a_1^*,a_2),
(a_2^*,a_3),\ldots,(a_{m-1}^*,a_m)$ are edges in $T^t(n,m)$. If we
delete the edges
$(b_1,b_1^*),(b_2,b_2^*),\ldots,(b_n,b_n^*);(b_1^*,b_2),(b_2^*,b_3),
\ldots,(b_{n-1}^*,b_n),(b_n^*,b_1)$ from $T^t(n,m)$, the triangular
lattice with cylindrical boundary condition, denoted by $T^c(n,m)$,
is obtained (see Figure 1(b)). If we delete the edges
$(a_1,a_1^*),(a_2,a_2^*),\ldots,(a_m,a_m^*), (a_1^*,a_2),\\
(a_2^*,a_3),\ldots,(a_{m-1}^*,a_m)$ from $T^c(n,m)$, the triangular
lattice with free boundary condition, denoted by $T^f(n,m)$ is
obtained (see Figure 1(c)). The asymptotic number of perfect
matchings of $T^t(n,m)$ can be found in Wu \cite{WFY06}. The
asymptotic number of spanning trees of $T^t(n,m)$ (resp. $T^c(n,m)$
and $T^f(n,m)$) was obtained by Shrock and Wu \cite{SW00} (resp. by
Yan and Zhang \cite{YZ08}).
\begin{theorem}
For the triangular lattices $T^t(n,m), T^c(n,m)$, and $T^f(n,m)$
with toroidal, cylindrical, and free boundary conditions,
$$\lim_{n,m\rightarrow \infty}\frac{\mathcal E(T^t(n,m))}{nm}=\lim_{n,m\rightarrow \infty}\frac{\mathcal E(T^c(n,m))}{nm}=
\lim_{n,m\rightarrow \infty}\frac{\mathcal
E(T^f(n,m))}{nm}$$$$=\frac{1}{2\pi^2}\int_0^{2\pi}\int_0^{2\pi}|\cos
x+\cos y+\cos(x+y)|{\rm d}_x{\rm d}_y(\approx 2.065),$$ that is, the
triangular lattices $T^t(n,m), T^c(n,m)$, and $T^f(n,m)$ with
toroidal, cylindrical, and free boundary conditions have the same
asymptotic energy ($\approx 2.065mn$).
\end{theorem}
\begin{proof}
By definitions of $T^t(n,m), T^c(n,m)$, and $T^f(n,m)$, $T^f(n,m)$
and $T^c(n,m)$ are spanning subgraph of $T^t(n,m)$. Moreover, almost
all vertices of $T^f(n,m)$ or $T^c(n,m)$ are $6$. Hence, by
Corollary 2.5,
$$\lim_{n,m\rightarrow \infty}\frac{\mathcal E(T^t(n,m))}{nm}=
\lim_{n,m\rightarrow \infty}\frac{\mathcal E(T^c(n,m))}{nm}=
\lim_{n,m\rightarrow \infty}\frac{\mathcal E(T^f(n,m))}{nm}.$$

It suffices to prove that
$$\lim_{n,m\rightarrow \infty}\frac{\mathcal E(T^t(n,m))}{nm}=\frac{1}{2\pi^2}\int_0^{2\pi}\int_0^{2\pi}|\cos
x+\cos y+\cos(x+y)|{\rm d}_x{\rm d}_y\approx 2.065.$$

Let $A(C_m)$ be the adjacency matrix of the cycle $C_m$. By a
suitable labelling of vertices of $T^t(n,m)$, the adjacency matrix
$A(T^t(n,m))$ of $T^t(n,m)$ has the following form:
$$A(T^t(n,m))=\left(
\begin{array}{cccccc}
A(C_m) & I_m+B_m & 0 & \cdots & 0 & I_m+B_m^T\\
I_m+B_m^T & A(C_m) & I_m+B_m & \cdots & 0 & 0 \\
0 & I_m+B_m^T & A(C_m) & \cdots & 0 & 0\\
\vdots & \vdots & \ddots & \ddots & \ddots & \vdots\\
0 & 0 & 0 & \cdots & A(C_m) & I_m+B_m\\
I_m+B_m & 0 & 0 & \cdots & I_m+B_m^T & A(C_m)
\end{array}
\right)_{n\times n}
$$
$$=I_n\otimes A(C_m)+B_n\otimes (I_m+B_m)+B_n^T\otimes (I_m+B_m^T),$$
where $I_n$ is the identity matrix of order $n$, $M\otimes N$
denotes the tensor product of two matrices $M$ and $N$, and
$$B_n=\left(
\begin{array}{cccccc}
0 & 1 & 0 & \cdots & 0 & 0\\
0 & 0 & 1 & \cdots & 0 & 0 \\
\vdots & \vdots & \ddots & \ddots & \ddots & \vdots\\
0 & 0 & 0 & \cdots & 0 & 1\\
1 & 0 & 0 & \cdots & 0 & 0
\end{array}
\right)_{n\times n}.$$ Note that $A(C_n)=B_n+B_n^T$. Hence
$$A(T^t(n,m))=I_n\otimes (B_m+B_m^T)+B_n\otimes (I_m+B_m)+B_n^T\otimes (I_m+B_m^T).$$
Let $\{1=g^0,g^1,\ldots,g^{n-1}\}$ be the cyclic group of order $n$.
Obviously, $\rho: g^i\rightarrow B_n^i$ for $0\leq i\leq n-1$ is a
representation of this group. Note that the cyclic group of order
$n$ has exactly $n$ (linear) characters $\chi_i$
($i=0,1,\ldots,n-1$) such that $\chi_i(g)=\omega_n^i$, where
$\omega_n$ is the $n$th root of unitary. Hence there exists an
invertible matrix $Q_n=(\frac{\omega_n^{ij}}{\sqrt n})_{0\leq
i,j\leq n-1}$ such that
$Q_n^{-1}B_nQ_n=diag(1,\omega_n,\ldots,\omega_n^{n-1})=: D_n$. Since
$B_n^T=B_n^{-1}$ and $Q_n^T=Q_n^{-1}$,
$Q_n^{-1}B_n^{T}Q_n=diag(1,\omega_n^{-1},\ldots,\omega_n^{-(n-1)})=:
D_n^{-1}$. Hence
$$(Q_n^{-1}\otimes Q_m^{-1})A(T^t(n,m))(Q_n\otimes Q_m)$$
$$=(Q_n^{-1}\otimes Q_m^{-1})[I_n\otimes (B_m+B_m^T)+B_n\otimes (I_m+B_m)+
B_n^T\otimes (I_m+B_m^T)](Q_n\otimes Q_m)$$
$$=I_n\otimes (D_m+D_m^{-1})+D_n\otimes (I_m+D_m)+D_n^{-1}\otimes (I_m+D_m^{-1}).$$
It is not difficult to see that $I_n\otimes
(D_m+D_m^{-1})+D_n\otimes (I_m+D_m)+D_n^{-1}\otimes (I_m+D_m^{-1})$
is a diagonal matrix whose diagonal entries are
$\omega_m^j+\omega_m^{-j}+\omega_n^i+\omega_n^{-i}+\omega_n^i\omega_m^j+\omega_n^{-i}\omega_m^{-j}
=2\cos\frac{2i\pi}{n}+2\cos\frac{2j\pi}{m}+2\cos(\frac{2i\pi}{n}+\frac{2j\pi}{m}),
0\leq i\leq n-1, 0\leq j\leq m-1$. This implies that the eigenvalues
of $A(T^t(n,m))$ are
$2\cos\frac{2i\pi}{n}+2\cos\frac{2j\pi}{m}+2\cos(\frac{2i\pi}{n}+\frac{2j\pi}{m}),
0\leq i\leq n-1, 0\leq j\leq m-1$. By the definition of the energy,
$$\mathcal E(T^t(n,m))=\sum_{i=0}^{n-1}\sum_{j=0}^{m-1}|2\cos\frac{2i\pi}{n}+
2\cos\frac{2j\pi}{m}+2\cos(\frac{2i\pi}{n}+\frac{2j\pi}{m})|.$$ So
$$\lim\limits_{n,m\rightarrow \infty}\frac{\mathcal
E(T^t(n,m))}{nm}= \lim_{n,m\rightarrow
\infty}\frac{1}{nm}\sum_{i=0}^{n-1}\sum_{j=0}^{m-1}|2\cos\frac{2i\pi}{n}+
2\cos\frac{2j\pi}{m}+2\cos(\frac{2i\pi}{n}+\frac{2j\pi}{m})|$$
$$=\frac{1}{2\pi^2}\int_0^{2\pi}\int_0^{2\pi}|\cos
x+\cos y+\cos(x+y)|{\rm d}_x{\rm d}_y\approx 2.065$$ and we complete
the proof of the theorem.
\end{proof}
\begin{figure}[htbp]
  \centering
 \scalebox{0.7}{\includegraphics{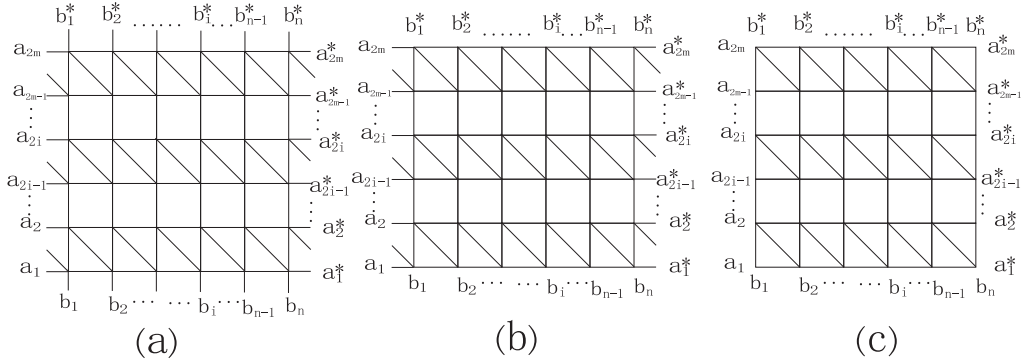}}
  \caption{\ (a)\ The $3^3.4^2$ lattice $S^t(n,2m)$ with toroidal boundary
  condition;
  \ (b)\ the $3^3.4^2$ lattice $S^c(n,2m)$ with cylindrical boundary condition;
  \ (c)\ the $3^3.4^2$ lattice $S^f(n,2m)$ with cylindrical boundary condition.}
\end{figure}

\subsection{The $3^3.4^2$ lattice}

The $3^3.4^2$ lattice $S^t(n,2m)$ with toroidal boundary condition
can be constructed by starting with a $2m\times n$ square lattice
and adding a diagonal edge connecting the vertices in, say, the
upper left to the lower right corners of each square in every other
row as shown in Figure 2(a), where
$a_1=b_1,a_{2m}=b_1^*,a_1^*=b_n,a_{2m}^*=b_n^*$, and
$(a_1,a_1^*),(a_2,a_2^*),\ldots,(a_{2m},a_{2m}^*),(b_1,b_1^*),(b_2,b_2^*),
\ldots,(b_n,b_n^*),(a_1,a_2^*),(a_3,a_4^*),\ldots,(a_{2m-1},a_{2m}^*)$
are edges in $S^t(n,2m)$. If we delete edges
$(b_1,b_1^*),(b_2,b_2^*),\ldots,(b_n,b_n^*)$ from $S^t(n,2m)$, the
$3^3.4^2$ lattice $S^c(n,2m)$ with cylindrical boundary condition is
obtained (see Figure 2(b)). If we delete edges
$(a_1,a_1^*),(a_2,a_2^*),\ldots,(a_{2m},a_{2m}^*),(a_1,a_2^*),(a_3,a_4^*),\ldots,(a_{2m-1},a_{2m}^*)$
from $S^c(n,2m)$, the $3^3.4^2$ lattice $S^f(n,2m)$ with free
boundary condition is obtained (see Figure 2(c)). The asymptotic
number of spanning trees of $S^t(n,2m)$ (resp. $S^c(n,2m)$ and
$T^f(n,2m)$) was obtained by Chang and Wang \cite{CW06} (resp. by
Yan and Zhang \cite{YZ08}).
\begin{theorem}
For the $3^3.4^2$ lattices $S^t(n,2m), S^c(n,2m)$, and $S^f(n,2m)$
with toroidal, cylindrical, and free boundary conditions,
$$\lim_{n,m\rightarrow \infty}\frac{\mathcal E(S^t(n,2m))}{2nm}=
\lim_{n,m\rightarrow \infty}\frac{\mathcal E(S^c(n,2m))}{2nm}=
\lim_{n,m\rightarrow \infty}\frac{\mathcal
E(S^f(n,2m))}{2nm}$$$$=\frac{1}{8\pi^2}\int_0^{2\pi}\int_0^{2\pi}(|2\cos
x+\sqrt{3+2\cos x+2\cos y+2\cos (x+y)}|+$$$$|2\cos x-\sqrt{3+2\cos
x+2\cos y+2\cos (x+y)}|){\rm d}_x{\rm d}_y\approx 1.8471,$$ that is,
the $3^3.4^2$ lattices $S^t(n,2m), S^c(n,2m)$, and $S^f(n,2m)$ with
toroidal, cylindrical, and free boundary conditions have the same
asymptotic energy ($\approx3.6942mn$).
\end{theorem}
\begin{proof}
By definitions of $S^t(n,2m), S^c(n,2m)$, and $S^f(n,2m)$,
$S^f(n,2m)$ and $S^c(n,2m)$ are spanning subgraphs of $S^t(n,2m)$.
Moreover, almost all vertices of $S^f(n,m)$ or $S^c(n,m)$ are of
degree $5$. Hence, by Corollary 2.5,
$$\lim_{n,m\rightarrow \infty}\frac{\mathcal E(S^t(n,m))}{2nm}=
\lim_{n,m\rightarrow \infty}\frac{\mathcal E(S^c(n,m))}{2nm}=
\lim_{n,m\rightarrow \infty}\frac{\mathcal E(S^f(n,2m))}{2nm}.$$

It suffices to prove that
$$\lim_{n,m\rightarrow \infty}\frac{\mathcal E(S^t(n,2m))}{2nm}=\frac{1}{8\pi^2}\int_0^{2\pi}\int_0^{2\pi}(|2\cos
x+\sqrt{3+2\cos x+2\cos y+2\cos (x+y)}|+$$$$|2\cos x-\sqrt{3+2\cos
x+2\cos y+2\cos (x+y)}|{\rm d}_x{\rm d}_y\approx 1.8471.$$

Let $A(C_{2m})$ be the adjacency matrix of the cycle $C_{2m}$. By a
suitable labelling of vertices of $S^t(n,2m)$, the adjacency matrix
$A(S^t(n,2m))$ of $S^t(n,2m)$ has the following form:
$$A(S^t(n,m))=\left(
\begin{array}{cccccc}
A(C_{2m}) & I_{2m}+F_{2m} & 0 & \cdots & 0 & I_{2m}+F_{2m}^T\\
I_{2m}+F_{2m}^T & A(C_{2m}) & I_{2m}+F_{2m} & \cdots & 0 & 0 \\
0 & I_{2m}+F_{2m}^T & A(C_{2m}) & \cdots & 0 & 0\\
\vdots & \vdots & \ddots & \ddots & \ddots & \vdots\\
0 & 0 & 0 & \cdots & A(C_{2m}) & I_{2m}+F_{2m}\\
I_{2m}+F_{2m} & 0 & 0 & \cdots & I_{2m}+F_{2m}^T & A(C_{2m})
\end{array}
\right)_{n\times n}
$$
$$=I_n\otimes A(C_{2m})+B_{n}\otimes (I_{2m}+F_{2m})+B_n^T\otimes (I_{2m}+F_{2m}^T),$$
where $I_n$ is the identity matrix of order $n$, $M\otimes N$
denotes the tensor product of two matrices $M$ and $N$,
$B_n=(b_{ij})_{n\times n}$ such that $b_{ij}=1$ if
$(i,j)=(1,2),(2,3),\ldots,(n-1,n),(n,1)$ and $b_{ij}=0$ otherwise,
and $F_{2m}=(f_{ij})_{2m\times 2m}$ such that $f_{ij}=1$ if
$(i,j)=(2,1),(4,3),(6,5),\ldots,(2m,2m-1)$ and $f_{ij}=0$ otherwise.

Note that $A(C_{2m})=I_m\otimes \left(
\begin{array}{cc}
0 & 1\\
1 & 0
\end{array}
\right)+B_{m}\otimes \left(
\begin{array}{cc}
0 & 0\\
1 & 0
\end{array}
\right)+B_{m}^T\otimes\left(
\begin{array}{cc}
0 & 1\\
0 & 0
\end{array}
\right)$. Hence
$$A(S^t(n,2m))=I_n\otimes I_m\otimes \left(
\begin{array}{cc}
0 & 1\\
1 & 0
\end{array}
\right)+I_n\otimes B_{m}\otimes \left(
\begin{array}{cc}
0 & 0\\
1 & 0
\end{array}
\right)+I_n\otimes B_{m}^T\otimes\left(
\begin{array}{cc}
0 & 1\\
0 & 0
\end{array}
\right)$$$$+B_n\otimes I_m\otimes \left(
\begin{array}{cc}
1 & 0\\
1 & 1
\end{array}
\right)+B_n^T\otimes I_m\otimes \left(
\begin{array}{cc}
1 & 1\\
0 & 1
\end{array}
\right).
$$
Let $Q_n=(\frac{\omega_n^{ij}}{\sqrt n})_{0\leq i,j\leq n-1}$. Then
$Q_n^{-1}B_nQ_n=diag(1,\omega_n,\ldots,\omega_n^{n-1})=: D_n$. Since
$B_n^T=B_n^{-1}$ and $Q_n^T=Q_n^{-1}$,
$Q_n^{-1}B_n^{T}Q_n=diag(1,\omega_n^{-1},\ldots,\omega_n^{-(n-1)})=:
D_n^{-1}$. Hence
$$(Q_n^{-1}\otimes Q_m^{-1}\otimes I_2)A(S^t(n,2m))(Q_n\otimes Q_m\otimes I_2)$$
$$=I_n\otimes I_m\otimes \left(
\begin{array}{cc}
0 & 1\\
1 & 0
\end{array}
\right)+I_n\otimes D_m\otimes \left(
\begin{array}{cc}
0 & 0\\
1 & 0
\end{array}
\right)+I_n\otimes D_m^{-1}\otimes \left(
\begin{array}{cc}
0 & 1\\
0 & 0
\end{array}
\right)$$$$+D_n\otimes I_m\otimes \left(
\begin{array}{cc}
1 & 0\\
1 & 1
\end{array}
\right)+D_n^{-1}\otimes I_m\otimes \left(
\begin{array}{cc}
1 & 1\\
0 & 1
\end{array}
\right).$$ It is not difficult to see that the above matrix is a
block diagonal matrix whose diagonal blocks are $\left(
\begin{array}{cc}
\omega_n^i+\omega_n^{-i} & 1+\omega_n^{-i}+\omega_m^j\\
1+\omega_n^{i}+\omega_m^{-j} & \omega_n^i+\omega_n^{-i}
\end{array}
\right), 0\leq i\leq n-1, 0\leq j\leq m-1$. Hence the eigenvalues of
$A(S^t(n,2m))$ are:
$$2\cos\frac{2i\pi}{n}\pm\sqrt{3+2\cos\frac{2i\pi}{n}+2\cos\frac{2j\pi}{m}+2\cos\left(\frac{2i\pi}{n}+\frac{2j\pi}{m}\right)},
0\leq i\leq n-1, 0\leq j\leq m-1.$$ Set
$$\lambda_1(i,j)=2\cos\frac{2i\pi}{n}+\sqrt{3+2\cos\frac{2i\pi}{n}+2\cos\frac{2j\pi}{m}+2\cos\left(\frac{2i\pi}{n}+\frac{2j\pi}{m}\right)},$$
$$\lambda_2(i,j)=2\cos\frac{2i\pi}{n}-\sqrt{3+2\cos\frac{2i\pi}{n}+2\cos\frac{2j\pi}{m}+2\cos\left(\frac{2i\pi}{n}+\frac{2j\pi}{m}\right)}.$$

By the definition of the energy, $$\mathcal
E(S^t(n,2m))=\sum\limits_{i=0}^{n-1}\sum\limits_{j=0}^{m-1}|(\lambda_1(i,j)|+|\lambda_2{i,j}|).$$
So $$\lim\limits_{n,m\rightarrow \infty}\frac{\mathcal
E(S^t(n,2m))}{2nm}=\frac{1}{8\pi^2}\int_0^{2\pi}\int_0^{2\pi}(|2\cos
x+\sqrt{3+2\cos x+2\cos y+2\cos (x+y)}|+$$$$|2\cos x-\sqrt{3+2\cos
x+2\cos y+2\cos (x+y)}|){\rm d}_x{\rm d}_y\approx 1.8471$$ and we
complete the proof of the theorem.
\end{proof}
\begin{figure}[htbp]
  \centering
 \scalebox{0.7}{\includegraphics{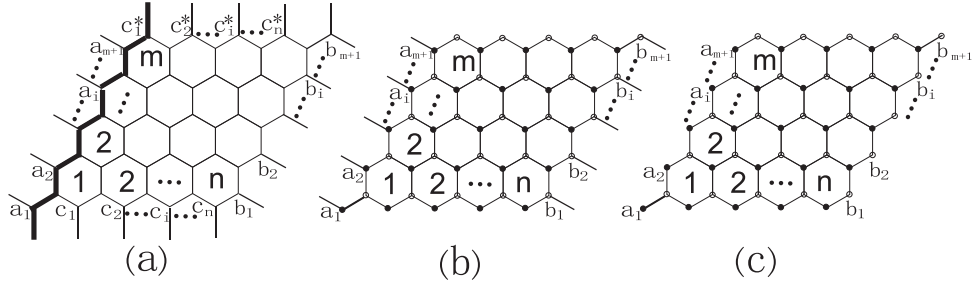}}
  \caption{\ (a)\ The hexagonal lattice $H^t(n,m)$ with toroidal boundary
  condition;
  \ (b)\ the hexagonal lattice $H^c(n,m)$ with cylindrical boundary condition;
  \ (c)\ the hexagonal lattice $H^f(n,m)$ with cylindrical boundary condition.}
\end{figure}
\subsection{The hexagonal lattice}
The hexagonal lattices with toroidal and cylindrical boundary
conditions, denoted by $H^t(n,m)$ and $H^c(n,m)$, are illustrated in
Figure 3(a) and Figure 3(b), where $(a_1,b_1), (a_2,b_2), \ldots,
(a_{m+1},b_{m+1}), (a_1,c^*_1), (c_1,c^*_2), (c_2,c^*_3)$, $\ldots,
(c_{n-1},c^*_n),\\ (c_n,b_{m+1})$ are edges in $H^t(n,m)$, and
$(a_1,b_1), (a_2,b_2), \ldots, (a_{m+1},b_{m+1})$ are edges in
$H^c(n,m)$. If we delete edges $(a_1,b_1), (a_2,b_2), \ldots,
(a_{m+1},b_{m+1})$ from $H^c(n,m)$, then the hexagonal lattice,
denoted by $H^f(n,m)$, with free boundary condition is obtained (see
Figure 3(c)). The asymptotic number of perfect matchings of
$H^t(n,m)$ can be found in Wu \cite{WFY06}. The asymptotic number of
spanning trees of $H^t(n,m)$ (resp. $H^c(n,m)$ and $H^f(n,m)$) was
obtained by Shrock and Wu \cite{SW00} (resp. by Yan and Zhang
\cite{YZ08}).
\begin{theorem}
For the hexagonal lattices $H^t(n,m), H^c(n,m)$, and $H^f(n,m)$ with
toroidal, cylindrical, and free boundary conditions,
$$\lim_{n,m\rightarrow \infty}\frac{\mathcal E(H^t(n,m))}{2nm}=
\lim_{n,m\rightarrow \infty}\frac{\mathcal E(H^c(n,m))}{2nm}=
\lim_{n,m\rightarrow \infty}\frac{\mathcal
E(H^f(n,m))}{2nm}$$$$=\frac{1}{4\pi^2}\int_0^{2\pi}\int_0^{2\pi}\sqrt{3+2\cos
x+2\cos y+2\cos (x+y)}{\rm d}_x{\rm d}_y\approx 1.5746,$$ that is,
the hexagonal lattices $H^t(n,m), H^c(n,m)$, and $H^f(n,m)$ with
toroidal, cylindrical, and free boundary conditions have the same
asymptotic energy ($\approx3.1492mn$).
\end{theorem}
\begin{proof}
By definitions of $H^t(n,m), H^c(n,m)$, and $H^f(n,m)$, $H^f(n,m)$
and $H^c(n,m)$ are spanning subgraphs of $H^t(n,m)$. Moreover,
almost all vertices of $H^f(n,m)$ or $H^c(n,m)$ are $3$. Hence, by
Corollary 2.5,
$$\lim_{n,m\rightarrow \infty}\frac{\mathcal E(H^t(n,m))}{2nm}=
\lim_{n,m\rightarrow \infty}\frac{\mathcal E(H^c(n,m))}{2nm}=
\lim_{n,m\rightarrow \infty}\frac{\mathcal E(H^f(n,m))}{2nm}.$$

It suffices to prove that
$$\lim_{n,m\rightarrow \infty}\frac{\mathcal E(H^t(n,m))}{2nm}
=\frac{1}{4\pi^2}\int_0^{2\pi}\int_0^{2\pi}\sqrt{3+2\cos x+2\cos
y+2\cos (x+y)}{\rm d}_x{\rm d}_y\approx 1.5746.$$

Let $A(C_{2m})$ be the adjacency matrix of the cycle $C_{2m}$. By a
suitable labelling of vertices of $H^t(n-1,m-1)$, the adjacency
matrix $A(H^t(n-1,m-1))$ of $H^t(n-1,m-1)$ has the following form:
$$A(H^t(n-1,m-1))=\left(
\begin{array}{cccccc}
A(C_{2m}) & F_{2m} & 0 & \cdots & 0 & F_{2m}^T\\
F_{2m}^T & A(C_{2m}) & F_{2m} & \cdots & 0 & 0 \\
0 & F_{2m}^T & A(C_{2m}) & \cdots & 0 & 0\\
\vdots & \vdots & \ddots & \ddots & \ddots & \vdots\\
0 & 0 & 0 & \cdots & A(C_{2m}) & F_{2m}\\
F_{2m} & 0 & 0 & \cdots & F_{2m}^T & A(C_{2m})
\end{array}
\right)_{n\times n}
$$
$$=I_n\otimes A(C_{2m})+B_{n}\otimes F_{2m}+B_n^T\otimes F_{2m}^T,$$
where $I_n$ is the identity matrix of order $n$, $M\otimes N$
denotes the tensor product of two matrices $M$ and $N$,
$B_n=(b_{ij})_{n\times n}$ such that $b_{ij}=1$ if
$(i,j)=(1,2),(2,3),\ldots,(n-1,n),(n,1)$ and $b_{ij}=0$ otherwise,
and $F_{2m}=(f_{ij})_{2m\times 2m}$ such that $f_{ij}=1$ if
$(i,j)=(2,1),(4,3),(6,5),\ldots,(2m,2m-1)$ and $f_{ij}=0$ otherwise.

Note that $A(C_{2m})=I_m\otimes \left(
\begin{array}{cc}
0 & 1\\
1 & 0
\end{array}
\right)+B_{m}\otimes \left(
\begin{array}{cc}
0 & 0\\
1 & 0
\end{array}
\right)+B_{m}^T\otimes\left(
\begin{array}{cc}
0 & 1\\
0 & 0
\end{array}
\right)$. Hence
$$A(H^t(n-1,m-1))=I_n\otimes I_m\otimes \left(
\begin{array}{cc}
0 & 1\\
1 & 0
\end{array}
\right)+I_n\otimes B_{m}\otimes \left(
\begin{array}{cc}
0 & 0\\
1 & 0
\end{array}
\right)+I_n\otimes B_{m}^T\otimes\left(
\begin{array}{cc}
0 & 1\\
0 & 0
\end{array}
\right)$$$$+B_n\otimes I_m\otimes \left(
\begin{array}{cc}
0 & 0\\
1 & 0
\end{array}
\right)+B_n^T\otimes I_m\otimes \left(
\begin{array}{cc}
0 & 1\\
0 & 0
\end{array}
\right).
$$
Let $Q_n=(\frac{\omega_n^{ij}}{\sqrt n})_{0\leq i,j\leq n-1}$. Then
$Q_n^{-1}B_nQ_n=diag(1,\omega_n,\ldots,\omega_n^{n-1})=: D_n$. Since
$B_n^T=B_n^{-1}$ and $Q_n^T=Q_n^{-1}$,
$Q_n^{-1}B_n^{T}Q_n=diag(1,\omega_n^{-1},\ldots,\omega_n^{-(n-1)})=:
D_n^{-1}$. Hence
$$(Q_n^{-1}\otimes Q_m^{-1}\otimes I_2)A(H^t(n-1,m-1))(Q_n\otimes Q_m\otimes I_2)$$
$$=I_n\otimes I_m\otimes \left(
\begin{array}{cc}
0 & 1\\
1 & 0
\end{array}
\right)+I_n\otimes D_m\otimes \left(
\begin{array}{cc}
0 & 0\\
1 & 0
\end{array}
\right)+I_n\otimes D_m^{-1}\otimes \left(
\begin{array}{cc}
0 & 1\\
0 & 0
\end{array}
\right)$$$$+D_n\otimes I_m\otimes \left(
\begin{array}{cc}
0 & 0\\
1 & 0
\end{array}
\right)+D_n^{-1}\otimes I_m\otimes \left(
\begin{array}{cc}
0 & 1\\
0 & 0
\end{array}
\right).$$ It is not difficult to see that the above matrix is a
block diagonal matrix whose diagonal blocks are $\left(
\begin{array}{cc}
0 & 1+\omega_n^{-i}+\omega_m^j\\
1+\omega_n^{i}+\omega_m^{-j} & 0
\end{array}
\right), 0\leq i\leq n-1, 0\leq j\leq m-1$. Hence the eigenvalues of
$A(H^t(n-1,m-1))$ are:
$\pm\sqrt{3+2\cos\frac{2i\pi}{n}+2\cos\frac{2j\pi}{m}+2\cos\left(\frac{2i\pi}{n}+\frac{2j\pi}{m}\right)},
0\leq i\leq n-1, 0\leq j\leq m-1.$

By the definition of the energy,
$$\mathcal
E(H^t(n,m))=2\sum\limits_{i=0}^{n}\sum\limits_{j=0}^{m}\sqrt{3+2\cos\frac{2i\pi}{n+1}+
2\cos\frac{2j\pi}{m+1}+2\cos\left(\frac{2i\pi}{n+1}+\frac{2j\pi}{m+1}\right)}.$$
So $$\lim\limits_{n,m\rightarrow \infty}\frac{\mathcal
E(H^t(n,m))}{2(n+1)(m+1)}=\frac{1}{4\pi^2}\int_0^{2\pi}\int_0^{2\pi}\sqrt{3+2\cos
x+2\cos y+2\cos (x+y)}{\rm d}_x{\rm d}_y\approx 1.5746$$ and we
complete the proof of the theorem.
\end{proof}

We would like to point out that the result of hexagonal lattice with
toroidal boundary condition has been obtained in \cite{H08} by a
different approach.
\begin{remark}
For the triangular, $3^3.4^2$, and hexagonal lattices, we have
considered the three boundary conditions:  the toroidal,
cylindrical, and free boundary conditions in Theorem 2.8-2.10,
respectively. By a similar idea, we can consider another two
boundary conditions: the Mobius-band and Klein-bottle boundary
conditions, and show that the triangular (resp. $3^3.4^2$ and
hexagonal) lattices with these five boundary conditions have the
same asymptotic energy.
\end{remark}


\section{CONCLUDING REMARKS}
In this paper, we showed that for many types of lattices the energy
per vertex of the plane lattices is independent of the boundary
conditions. It is no difficulty to see that the conclusion is true
in general. In fact our approach can be used widely. By using this
conclusion we can convert some harder problem to easy one and get
some results simultaneously. For example dealing with the problem of
the asymptotic energy of the hexagonal lattice with the free
boundary is not an easy task but we deduced it in a simple way. On
the other hand, for the entropy of dimers the result is not true. In
fact, Yan, Yeh, and Zhang \cite{YYZ08} showed that for the dimer
problem the hexagonal lattices with cylindrical and toroidal
boundary have different entropies.


\vskip1cm \noindent
\newcounter{cankao}
\begin{list}
{[\arabic{cankao}]}{\usecounter{cankao}\itemsep=0cm} \centerline{\bf
References}

\bibitem{CW98}
N. J. Calkin and H. S. Wilf, \textit{The number of independent sets
in a grid graph}, SIAM J. DISCRETE MATH. Vol. 11 (1998) 54¨C-60.

\bibitem{CW06}
S.-C. Chang and W.Wang, \textit{Spanning trees on lattices and
integral identities}, J. Phys. A: Math. Gen., 39 (2006)
10263--10275.

\bibitem{DS07}
J. Day, W. So, \textit{Singular value inequality and graph energy
change}, Electronic Journal of Linear Algebra, 16 (2007), 291--299.

\bibitem{DS08}
J. Day, W. So, \textit{Graph energy change due to edge deletion},
Linear Algebra and its Applications, 428 (2008), 2070--2078.

\bibitem{Fisher61}
M. E. Fisher, \textit{Statistical mechanics of dimers on a plane
lattice}, Phys. Rev., 124(1961), 1664--1672.

\bibitem{Gutm78}
I. Gutman, \textit{The energy of a graph}, Ber. Math. -Statist.
Sekt. Forschungszentrum Graz, 103 (1978), 1--22.

\bibitem{Gutm99}
I. Gutman, The energy of a graph: old and new results, Algebraic
Combinatorics and Applications (G¡§ossweinstein, 1999), Springer,
Berlin, 2001, pp. 196--211.




\bibitem{H08}
J. H. Hua, \textit{The spectra of some lattice graphs on surfaces},
Dissertation of Lanzhou University, 2008.

\bibitem{Kast61}
P. W. Kasteleyn, \textit{The statistics of dimers on a lattice I:
The number of dimer arrangements on a quadratic lattice}, Physica,
27(1961), 1209--1225.

\bibitem{Kast63}
P. W. Kasteleyn, \textit{Dimer statistics and phase transitions}, J.
Math. Phys., 4(1963), 287--293.

\bibitem{Kast67}
P. W. Kasteleyn, Graph Theory and Crystal Physics, Graph Theory and
Theoretical Physics (F.Harary, ed.), Academic Press, 1967, 43--110.

\bibitem{Kool01}
J. H. Koolen, V. Moulton, \textit{Maximal energy graphs}, Adv. Appl.
Math. 26 (2001), 47--52.







\bibitem{TF61}
H. N. V. Temperley, M. E. Fisher, \textit{Dimer problem in
statistical mechanics---an exact result}, Philosophical Magazine,
6(1961), 1061--1063.



\bibitem{SW00}
R. Shrock, F. Y. Wu, \textit{Spanning trees on graphs and lattices
in $d$ dimensions}, J. Phys. A: Math. Gen., 33 (2000), 3881--3902.


\bibitem{WFY06}
F. Y. Wu, \textit{Dimers on two-dimensional lattices}, Intern. J.
Modern Phys. B, 32 (2006), 5357--5371.

\bibitem{YYZ08}
W. G. Yan, Y.-N. Yeh, and F. J. Zhang, \textit{Dimer problem on the
cylinder and torus}, Physica A: Statistical Mechanics and its
Applications, 387 (2008), 6069--6078.

\bibitem{YZ08}
W. G. Yan, F. J. Zhang, \textit{Enumeration of spanning trees of
some symmetric graphs}, Preprint.

\end{list}
\end{document}